\documentclass[11pt]{article}%
\usepackage[dvips]{epsfig}
\usepackage{amsmath}
\usepackage{amsfonts}
\usepackage{amssymb}
\usepackage{graphicx}
\usepackage{amsmath,amsfonts,amssymb,mathrsfs,amscd}
\usepackage{psfrag}
\usepackage{bbm}
\usepackage{float}
\usepackage{latexsym}%
\newtheorem{theorem}{Theorem}[section]

\newtheorem{lemma}{Lemma}[section]

\numberwithin{equation}{section}

\begin{document}

\title{{\LARGE On  conformally flat cubic metrics with weakly isotropic scalar curvature}}
\author{{ Cuiling Ma and Xiaoling Zhang$^*$  \ \ } }
\maketitle

\begin{abstract}
  The conformal properties of metrics are meaningful in Riemannian and Finsler geometry, and cubic metrics are useful in physics and biology. In this paper, we study the conformally flat cubic metrics with weakly isotropic scalar curvature. We also prove that such metrics must be Minkowski metrics.
\vspace{-3mm}
\end{abstract}

%-------------------  First Head  -----------------------------------------

%===================Text=============================================

%\vspace{8true mm}

%×÷Õß
%\date{}
%\renewcommand{\baselinestretch}{1.5}\baselineskip 12pt
%\centerline{\footnotesize\rm (Department of Mathematics, Xinjiang
%Normal University, Urumqi {\rm 830054}, China)}

\footnotetext{\baselineskip 10pt
%$^\dag$ Corresponding author\\
This work was supported partially by National Natural Science Foundation of China(No.11961061,11761069).  \\
*corresponding author} \thispagestyle{empty}

\vspace{1mm}{\footnotesize \textbf{Keywords}:\hspace{2mm}   weakly isotropic scalar curvature, conformally flat, cubic metrics,
%¹Ø¼ü´Ê
}

%\no{\footnotesize{\bf MSC(2000):\hspace{2mm} 20M31, 33Q21}
\vspace{2mm} \baselineskip 15pt \renewcommand{\baselinestretch}{1.22}
\parindent=10.8pt

\section{ Introduction}

The first to treat the conformal theory of Finsler metrics generally was, to the best of authors' knowledge, M.S.Knebelman. He founded the condition, that the length of an arbitrary vector in the one metric is proportional to the length in the other metric, in turn implies that the difference of the two metrics is only a function of position \cite{Knebelman1929}. We shall call this result Knebelman's theorem. Later, Hashiguchi  \cite{Hashiguchi1976} gave a geometrical definition. He stated that conformation means the angles between any vector $V$ and a supporting element $y$ are equal. In Finsler geometry, the Weyl theorem \cite{Rund1959} states that the projective and conformal properties of a Finsler space determine the metric properties uniquely. So the study of conformal properties of a Finsler metric becomes more important and it has been a recent popular trend in Finsler geometry.

It is one hot issue that how to characterize conformally flat metrics in Finsler geometry. There are many important results in conformal Finsler geometry. Kang proved that any conformally flat Randers metric of scalar flag curvature is projectively flat. Moreover, such metrics are completely classified(see \cite{Kang2011}).
Chen-He-Shen \cite{Chen2015} have proved that conformally flat $(\alpha,\beta)$-metrics with constant flag curvature must be Riemannian or locally Minkowskian.
In \cite{Tayebi2019}, Tayebi-Razgordani-Najafi showed that the conformally flat cubic metric with relatively isotropic mean Landsberg curvature on a manifold $M$ of dimension $n\geq3$ is either a Riemannian metric or a locally Minkowski metric.

The theory of $m$-th root metrics has been developed by Shimada \cite{Shimada1979}, and applied to Biology as ecological metrics by Antonelli \cite{Antonelli1993}. Later, many scholars studied these metrics (\cite{Aranasiu2009},\cite{Brinzei2008},etc).
In \cite{Wegener1935}, Wegener studied cubic Finsler metrics of dimensions two and three. Wegener's paper is only an abstract of his PhD thesis without almost all calculations. Further, in \cite{Matsumoto1996}, Matsumoto gave an improved version of Wegener's paper.
Zhang-Xia \cite{Zhang2014} studied Einstein Matsumoto metrics. They obtained some equivalent conditions of $F$ and proved that non-Riemannian Matsumoto metric on a 3-dimensional manifold $M$ is an Einstein metric if and only if it is of zero flag curvature.
Yu-You \cite{Yu2010} showed that the spray coefficients $G^{i}$ and Ricci curvature $\textbf{Ric}$ of $m$-th root metric must be rational functions in $y$. Based on the result in \cite{Cheng2012}, Chen-Xia \cite{Chen2019} gave the explicit expressions of $Ricci$ curvature and scalar curvature of $(\alpha,\beta)$-metric which is conformally flat. They studied the conformally flat $(\alpha,\beta)$-metrics with weakly isotropic scalar curvature. And they showed that if $\phi(s)$ in $F$ is a polynomial of degree $n\geq2$, then scalar curvature vanishes.

In this paper, we mainly focus on the conformally flat cubic metrics with weakly isotropic scalar curvature and obtain the following results.

\begin{theorem}\label{0}\quad
 Let the $m$-th root metric $F$ be of weakly isotropic scalar curvature. Then the scalar curvature $\textbf{r}$ must vanish.
\end{theorem}

\begin{theorem}\label{1}\quad
Let $F=F(x,y)$ be a conformally flat cubic Finsler metric on a manifold $M$ of dimension $n\geq 3$. Suppose $F$ is of weakly isotropic scalar curvature, then $F$ must be locally Minkowskian.
\end{theorem}

%#############################################################################################

\section{ Preliminaries}
In this section, we mainly introduce several geometric quantities in Finsler geometry and several results which will be used later.

Let $M$ be an $n$-dimensional smooth manifold with $n \geq3$. The points in the tangent bundle $TM$ are denoted by $(x,y)$, where $x \in M$ and $y \in T_{x}M$. Let $(x^{i},y^{i})$ be the local coordinates of $TM$ with $y=y^{i} \frac{\partial}{\partial x^{i}}$. A Finsler metric on $M$ is a function $F: TM \longrightarrow [0,+\infty)$ such that\\
(1) $F$ is smooth in $TM \backslash \{0\}$;\\
(2) $F(x,\lambda y)=\lambda F(x,y)$ for any $\lambda \textgreater 0$;\\
(3) The fundamental quadratic form  $g=g_{ij}(x,y) dx^{i} \otimes dx^{j}$, where
$$g_{ij}(x,y)=[\frac{1}{2} F^2(x,y)]_{y^{i}y^{j}} $$
is positively definite. We use the notations: $F_{y^i}:=\frac{\partial F}{\partial y^i}$, $F_{x^i}:=\frac{\partial F}{\partial x^i}, F^{2}_{y^iy^j}:=\frac{\partial^{2} F^{2}}{\partial y^i \partial y^j}.$

Let $F$ be a Finsler metric on an $n$-dimensional manifold $M$ and $G^{i}$ be the geodesic coefficients of $F$, which are defined by
\begin{equation*}\label{G}
G^{i}:=\frac{1}{4} g^{ij} (F^{2}_{x^k y^j} y^k-F^{2}_{x^j} ),
\end{equation*}
where $(g^{ij})=(g_{ij})^{-1}.$ For any $x \in M$ and $y \in T_x M \backslash\{0\}$, the Riemann curvature $R_{y}:=R^{i}_{\,\,\,\,k}(x,y) \frac{\partial}{\partial x^{i}} \otimes dx^{k}$ is defined by \\
$$R^{i}_{\,\,\,\,k}:=2 G^{i}_{x^{k}}-G^{i}_{x^{j}y^{k}} y^{j}+2 G^{j} G^{i}_{y^{j}y^{k}}-G^{i}_{y^{j}} G^{j}_{y^{k}}.$$
The Ricci curvature $\textbf{Ric}$ is the trace of the Riemann curvature defined by
\begin{equation*}
\textbf{Ric}:=R^{k}_{\,\,\,\,k}.
\end{equation*}
The Ricci tensor is
$$\textbf{Ric}_{ij}:=\frac{1}{2} \textbf{Ric}_{y^{i}y^{j}}.$$
By the homogeneity of $\textbf{Ric}$, we have $\textbf{Ric}=\textbf{Ric}_{ij} y^{i}y^{j}.$
The scalar curvature $\textbf{r}$ of $F$ is defined as
\begin{equation}\label{2.1}
\textbf{r}:=g^{ij} \textbf{Ric}_{ij}.
\end{equation}
A Finsler metric is said to be of weakly isotropic scalar curvature if there exists a 1-form $\theta=\theta_{i}(x) y^{i}$ and a scalar function $\mu(x)$ such that
\begin{equation}\label{2.2}
\textbf{r}=n (n-1) (\frac{\theta}{F}+\mu).
\end{equation}

An $(\alpha,\beta)$-metric is a Finsler metric of the form
\begin{equation*}\label{2.3}
F=\alpha \phi(s),
\end{equation*}
where $\alpha=\sqrt{a_{ij}(x) y^i y^j}$ is a Riemannian metric, $\beta=b_{i}(x) y^{i}$ is a 1-form, $s:=\frac{\beta}{\alpha}$, $b:=\parallel \beta \parallel_{\alpha} < b_0$, $x \in M$. It has been proved that $F=\alpha \phi(s)$ is a positive definite Finsler metric if and only if $\phi=\phi(s)$ is a positive $C^{\infty}$ function on $(-b_0,b_0)$ and satisfying the following condition
\begin{equation}\label{2.4}
\phi(s)-s \phi'(s)+(B-s^2) \phi''(s)>0,\quad |s| \le b <b_0,
\end{equation}
where $B:=b^2$.

Let $F=\sqrt[3]{a_{ijk}(x) y^{i} y^{j} y^{k}}$ be a cubic metric on a manifold $M$ of dimension $n\geq 3.$ By choosing suitable non-degenerate quadratic form $\alpha=\sqrt{a_{ij}(x) y^{i} y^{j}}$ and one-form $\beta=b_{i}(x) y^{i}$, it can be written in the form
\begin{equation*}\label{2.5}
F=\sqrt[3]{p \beta \alpha^{2}+q \beta^{3}},
\end{equation*}
where $p$ and $q$ are real constants such that $p+q B\neq 0$ (see \cite{Matsumoto1979}). The above equation can be rewritten as
\begin{equation*}\label{2.6}
F=\alpha (p s+q s^{3})^{\frac{1}{3}},
\end{equation*}
which means that $F$ is also an $(\alpha,\beta)$-metric with $\phi(s)=(p s+q s^{3})^{\frac{1}{3}}$. Then, by (\ref{2.4}), we obtain
\begin{equation}\label{2.7}
\phi^{-5} [-p^2 B +p (4 p+3 q B) s^2]>0.
\end{equation}

Two Finsler metrics $F$ and $\widetilde{F}$ on a manifold $M$ are said to be conformally related if there is a scalar function $\kappa(x)$ on $M$ such that $F=e^{\kappa(x)}\widetilde{F}.$ Particularly, an $(\alpha, \beta)$-metric $F=\alpha \phi(\frac{\beta}{\alpha})$ is said to be conformally related to a Finsler metric $\widetilde{F}$ if $F=e^{\kappa(x)} \widetilde{F}$ with $\widetilde{F}=\widetilde{\alpha} \phi(\widetilde{s})=\widetilde{\alpha} \phi(\frac{\widetilde{\beta}}{\widetilde{\alpha}}).$ In the following, we always use symbols with tilde to denote the corresponding quantities of the metric $\widetilde{F}.$ Note that $\alpha=e^{\kappa}\widetilde{\alpha}$, $\beta=e^{\kappa}\widetilde{\beta},$ thus $\widetilde{s}=s.$

A Finsler metric which is conformally related to a locally Minkowski metric is said to be conformally flat. Thus, a conformally flat $(\alpha,\beta)$-metric $F$ has the form $F=e^{\kappa(x)} \widetilde{F}$, where $\widetilde{F}=\widetilde{\alpha} \phi(\frac{\widetilde{\beta}}{\widetilde{\alpha}})$ is a locally Minkowski metric.

Denoting
\begin{equation*}\begin{aligned}\label{rs}
&r_{ij}:=\frac{1}{2} (b_{i|j}+b_{j|i}),\quad s_{ij}:=\frac{1}{2} (b_{i|j}-b_{j|i}),\\
&r^{i}_{\,\,\,\,j}:=a^{il} r_{lj},\quad s^{i}_{\,\,\,\,j}:=a^{il} s_{lj},\\
&r_{j}:=b^{i} r_{ij},\quad r:=b^i r_i,\quad s_{j}:=b^{i} s_{ij},\\
& r_{00}:=r_{ij} y^{i} y^{j},\quad s^{i}_{\,\,\,\,0}:=s^{i}_{\,\,\,\,j} y^{j},\quad s_{0}:=s_{i} y^{i},\\
\end{aligned}\end{equation*}
where $b^{i}:=a^{ij} b_{j}$, $b_{i|j}$ denotes the covariant defferantiation with respect to $\alpha$.

Let $G^{i}$ and $G^{i}_{\alpha}$ denote the geodesic coefficients of $F$ and $\alpha$, respectively. The geodesic coefficients $G^{i}$ of $F=\alpha \phi(\frac{\beta}{\alpha})$ are related to $G_{\alpha}^{i}$ by
\begin{equation*}
G^{i}=G^{i}_{\alpha}+\alpha Q s^{i}_{\,\,\,\,0}+(-2 Q \alpha s_{0}+r_{00}) (\Psi b^{i}+\Theta \alpha^{-1} y^{i}),
\end{equation*}
where
\begin{equation*}\begin{aligned}\label{3.1}
&Q:=\frac{\phi^{'}}{\phi-s \phi^{'}},\quad
 \Theta:=\frac{\phi \phi^{'}-s (\phi \phi^{''}+\phi^{'} \phi^{'})}{2 \phi [(\phi-s \phi)+(B-s^2) \phi^{''}]},\\
  & \Psi:=\frac{\phi^{''}}{2 [(\phi-s \phi^{'})+(B-s^2) \phi^{''}]}.\\
\end{aligned}\end{equation*}

Assume that $F=\alpha \phi(\frac{\beta}{\alpha})$ is conformally related to a Finsler metric $\widetilde{F}=\widetilde{\alpha} \phi(\frac{\widetilde{\beta}}{\widetilde{\alpha}})$ on $M$, i.e. $F=e^{\kappa(x)} \widetilde{F}$, where $\kappa(x)$ is a scalar function on $M$. Then
\begin{equation*}\begin{aligned}\label{3.2}
&a_{ij}=e^{2k} \widetilde{a}_{ij}, b_{i}=e^{k} \widetilde{b}_{i},  \widetilde{b}:=\parallel \widetilde{\beta} \parallel_{\widetilde{\alpha}}=\sqrt{\widetilde{a}_{ij} \widetilde{b}^{i} \widetilde{b}^{j}}=b.
\end{aligned}\end{equation*}
Further, we have
\begin{equation*}\begin{aligned}\label{3.3}
& b_{i|j}=e^{k(x)} (\widetilde{b}_{i \parallel j}-\widetilde{b}_{j} k_{i}+\widetilde{b}_{l} k^{l} \widetilde{a}_{ij}),\\
&^{\alpha}\Gamma^{l}_{ij}=^{\widetilde{\alpha}}\widetilde{\Gamma}^l_{ij}+k_j \delta^l_i+k_i \delta^l_j-k^l \widetilde{a}_{ij},\\
&r_{ij}=e^{k(x)} \widetilde{r}_{ij}+\frac{1}{2} e^{k(x)} (-\widetilde{b}_{j} k_{i}-\widetilde{b}_{i} k_{j}+2 \widetilde{b}_{l} k^{l} \widetilde{a}_{ij}),\\
&s_{ij}=e^{k(x)} \widetilde{s}_{ij}+\frac{1}{2} e^{k(x)} (\widetilde{b}_{i} k_{j}-\widetilde{b}_{j} k_{i}),\\
&r_{i}=\widetilde{r}_{i}+\frac{1}{2} (\widetilde{b}_{l} k^{l} \widetilde{b}_{i}-b^2 k_{i}),r=e^{-k(x)} \widetilde{r},\\
&s_{i}=\widetilde{s}_{i}+\frac{1}{2} (b^2 k_{i}-\widetilde{b}_{l} k^{l} \widetilde{b}_{i}),\\
&r^{i}_{\,\,\,\,\,\,i}=e^{-k(x)} \widetilde{r}^{i}_{\,\,\,\,\,\,i}+e^{-k(x)} (n-1)\widetilde{b}_{i} k^{i},\\
&s^{j}_{\,\,\,\,\,\,i}=e^{-k(x)} \widetilde{s}^j_{\,\,\,\,\,\,i}+\frac{1}{2} e^{-k(x)} (\widetilde{b}^{j} k_{i}-\widetilde{b}_{i} k^{j}).
\end{aligned}\end{equation*}
Here $\widetilde{b}_{i \parallel j}$ denote the covariant derivatives of $\widetilde{b}_{i}$ with respect to $\widetilde{\alpha}$, $^{\alpha}\Gamma^{m}_{ij}$ and $^{\widetilde{\alpha}}\widetilde{\Gamma}^{m}_{ij}$ denote Levi-Civita connection with respect to $\alpha$ and $\widetilde{\alpha}$, respectively. In the following we adopt $k_{i}:=\frac{\partial k}{\partial x^{i}}$, $k_{ij}:=\frac{\partial^2k}{\partial x^{i}\partial x^{j}}$, $k^{i}:=\widetilde{a}^{ij} k_{j}$, $\widetilde{b}^{i}:=\widetilde{a}^{ij} \widetilde{b}_{j}$, $f:=\widetilde{b}_{i} k^{i}$, $f_{1}:=k_{ij} \widetilde{b}^{i} y^{j}$, $f_{2}:=k_{ij} \widetilde{b}^{i} \widetilde{b}^{j}$, $k_{0}:=k_{i} y^{i}$, $k_{00}:=k_{ij} y^{i} y^{j}$ and $\parallel\triangledown k\parallel^2:=\widetilde{a}^{ij} k_{i} k_{j}$.

\begin{lemma}\label{2}(\cite{Chen2019})\quad
 Let $F=e^{\kappa(x)} \widetilde{F}$, where $\widetilde{F}=\widetilde{\alpha} \phi(\frac{\widetilde{\beta}}{\widetilde{\alpha}})$ is locally $Minkowskian$. Then the $Ricci$ curvature of $F$ is determined by
\begin{equation*}\label{RIC}
\textbf{Ric}=D_{1} \parallel\triangledown k\parallel^2 \widetilde{\alpha}+D_{2} k_{0}^{2}+D_{3} k_{0} f \widetilde{\alpha}+D_{4} f^{2} \widetilde{\alpha}^{2}+D_{5} f_{1} \widetilde{\alpha}+D_{6} \widetilde{\alpha}^{2}+D_{7} k_{00},
\end{equation*}
where $D_{k}(k=1,...,7)$ are listed in Lemma 3.2 in \cite{Chen2019}.
\end{lemma}

\begin{lemma}\label{3}(\cite{Chen2019})\quad
 Let $F=e^{\kappa(x)} \widetilde{F}$, where $\widetilde{F}=\widetilde{\alpha} \phi (\frac{\widetilde{\beta}}{\widetilde{\alpha}})$ is locally Minkowskian. Then the scalar curvature of $F$ is determined by
\begin{equation*}\label{3.7}
\textbf{r}=\frac{1}{2} e^{-2\kappa} \rho^{-1} [\Sigma_{1}-(\tau+\eta \lambda^{2}) \Sigma_{2}-\frac{\lambda \eta}{\widetilde{\alpha}} \Sigma_{3}-\frac{\eta}{\widetilde{\alpha}^{2}} \Sigma_{4}],
\end{equation*}
where
\begin{equation*}\begin{aligned}\label{3.8}
&\tau:=\frac{\delta}{1+\delta B},\quad \eta:=\frac{\mu}{1+Y^{2}\mu},\quad  \lambda:=\frac{\varepsilon-\delta s}{1+\delta B},\\
&\delta :=\frac{\rho_{0}-\varepsilon^{2} \rho_{2}}{\rho},\quad \varepsilon :=\frac{\rho_{1}}{\rho_{2}},\quad \mu :=\frac{\rho_{2}}{\rho},\\
&Y :=\sqrt{A_{ij} Y^{i} Y^{j}},\quad A_{ij} :=a_{ij}+\delta b_{i} b_{j},\\
&\rho :=\phi(\phi-s \phi'),\quad \rho_{0} :=\phi\phi''+\phi'\phi',\\
&\rho_{1} :=-s (\phi\phi''+\phi'\phi')+\phi\phi',\quad \rho_{2} :=s[s(\phi\phi''+\phi'\phi')-\phi\phi'],\\
\end{aligned}\end{equation*}
and $\Sigma_{i}(i=1,...,4)$ are listed in the proof of Lemma 3.3 in \cite{Chen2019}.
\end{lemma}

\begin{lemma}\label{4}(\cite{Yu2010})\quad
Let $m$-th root metric in the form $F=\sqrt[m]{a_{i_1 i_2 \cdots i_m}(x) y^{i_1} y^{i_2} \cdots y^{i_m}}$ be a Finsler metric on a manifold of dimension $n$. Then the Ricci curvature of $F$ is a rational function in $y$.
\end{lemma}

%#########################################################################################################################################

\section{Proof of the main theorems}
In this section, we will prove the main theorems. Firstly, we give the proof of Theorem \ref{0}.

{\bf{The proof of Theorem \ref{0}}  }

For an $m$-th root metric $F=\sqrt[m]{a_{i_1 i_2 \cdots i_m}(x) y^{i_1} y^{i_2} \cdots y^{i_m}}$ on a manifold $M$, the inverse of the fundamental tensor of $F$ is given by (see \cite{Yu2010})
\begin{equation}\label{4.1}
g^{ij}=\frac{F^{m-2}}{m-1} A^{ij}+\frac{m-2}{m-1} \frac{y^i y^j}{F^2},
\end{equation}
where $A_{ij}=\frac{1}{m (m-1)} \frac{\partial^2 F^m}{\partial y^i \partial y^j}$, $(A^{ij})=(A_{ij})^{-1}$. Thus $F^2 g^{ij}$ are rational functions in $y$.

By Lemma \ref{4}, the Ricci curvature $\textbf{Ric}$ of $m$-th root metric is a rational function in $y$, then $\textbf{Ric}_{ij}:=\textbf{Ric}_{y^{i}y^{j}}$ are rational functions. According to (\ref{2.1}), we have
\begin{equation}\begin{aligned}\label{4.5}
F^2 \textbf{r} =F^2 g^{ij} \textbf{Ric}_{ij}.
\end{aligned}\end{equation}
This means that $F^2 \textbf{r}$ is a rational function in $y$.

On the other hand, if $F$ is of weakly isotropic scalar curvature, according to (\ref{2.2}), we obtain
\begin{equation*}\begin{aligned}\label{4.4}
F^2 \textbf{r} =n (n-1) \rho (\theta F+\mu F^2),
\end{aligned}\end{equation*}
where $F=\alpha (p s+q s^{3})^{\frac{1}{3}}$, $\theta$ is a 1-form and $\mu$ is a scalar function. The right of the above equation is a irrational function in $y$. Comparing it with (\ref{4.5}), we have $\textbf{r}=0$.
$\hfill\blacksquare$

{\bf{The proof of Theorem \ref{1}}  }

Assume that the conformally flat cubic metric $F$ is of weakly isotropic scalar curvature. Then, by Lemma \ref{3} and Theorem \ref{0}, we obtain $\textbf{r}=0$, i.e.,
\begin{equation*}\label{4.6}
\Sigma_{1}-(\tau+\eta \lambda^{2}) \Sigma_{2}-\frac{\lambda \eta}{\widetilde{\alpha}} \Sigma_{3}-\frac{\eta}{\widetilde{\alpha}^{2} } \Sigma_{4}=0.
\end{equation*}
Further, by detailed expressions of $\Sigma_{i} (i=1,\cdots 4)$, the above equation can be rewritten as
\begin{equation}\label{4.3}
\frac{B (4 p+3 q B) \kappa_{0}^{2}-4 (4 p+3 q B) \widetilde{\beta} \kappa_{0}  f+4 p \widetilde{\alpha}^{2} f^{2}}{(4 p+3 q B)^8 \widetilde{\alpha}^2 s^2 \gamma^7}+\frac{T}{\gamma^6}=0,
\end{equation}
where $\gamma:=p B \widetilde{\alpha}^{2}-(4 p+3 q B) \widetilde{\beta}^{2}$ and $T$ has no $\gamma^{-1}$.

Thus the first term of (\ref{4.3}) can be divided by $\gamma$. It means that there is a function $h(x)$ on $M$ such that
\begin{equation*}\label{B1}
B (4 p+3 q B) \kappa_{0}^{2}-4 (4 p+3 q B) \widetilde{\beta} \kappa_{0}  f+4 p \widetilde{\alpha}^{2} f^{2}=h(x) \gamma.
\end{equation*}
The above equation can be rewritten as
\begin{equation}\label{B2}
B (4 p+3 q B) \kappa_{0}^{2}-4 (4 p+3 q B) \widetilde{\beta} \kappa_{0} f +4 p  \widetilde{\alpha}^{2} f^{2}=h(x) [ p B \widetilde{\alpha}^{2}-(4 p+3 q B) \widetilde{\beta}^{2}].
\end{equation}
Differentiating (\ref{B2}) with $y^{i}$ yields
\begin{equation}\label{4.8}
B (4 p+3 q B) \kappa_{0} \kappa_{i}-2 (4 p+3 q B) (\widetilde{b}_{i} \kappa_{0}+ \widetilde{\beta} \kappa_{i}) f+4 p \widetilde{a}_{il} y^{l} f^{2}= h(x) [ p B\widetilde{a}_{il} y^{l}-(4 p+3 q B)  \widetilde{\beta} \widetilde{b}_{i}].
\end{equation}
Differentiating (\ref{4.8}) with $y^{j}$ yields
\begin{equation*}\label{4.9}
B (4 p+3 q B)\kappa_{i} \kappa_{j}-2 (4 p+3 q B) (\widetilde{b}_{i} \kappa_{j}+\widetilde{b}_{j} \kappa_{i}) f+4 p \widetilde{a}_{ij} f^{2}=h(x) [ p B\widetilde{a}_{ij}-(4 p+3 q B) \widetilde{b}_{i} \widetilde{b}_{j}].
\end{equation*}
Contracting the above with $\widetilde{b}^{i} \widetilde{b}^{j}$ yields
\begin{equation*}\label{4.10}
B f^{2} (8 p+9 q B)=3 B^{2} h(x) [p+q B].
\end{equation*}
Thus we have
\begin{equation}\label{4.11}
h(x)=\frac{(8 p+9 q B) f^{2}}{3 B (p+q B)}.
\end{equation}

Substituting (\ref{4.11}) into (\ref{4.8}) and contracting (\ref{4.8}) with $\widetilde{b}^{i}$ yield
\begin{equation}\label{4.12}
(4 p+3 q B)f ( f\widetilde{\beta}-B \kappa_{0} )=0.
\end{equation}
Furthermore, by (\ref{2.7}) and $4 p+3 q B\neq 0$, we have $f ( f\widetilde{\beta}-B \kappa_{0} )=0$.\\

Case I: $f=0$. It means $h(x)=0$ by (\ref{4.11}). Thus, one have that $\kappa_{i}=0$ by (\ref{B2}), which means
\begin{equation*}\label{4.14}
\kappa=constant.
\end{equation*}

Case II: $f\neq 0$. It implies that $ f\widetilde{\beta}-B \kappa_{0}=0$. Substituting it into (\ref{B2}), we obtain
\begin{equation*}\label{4.13}
\widetilde{\beta}^2=-B \widetilde{\alpha}^2,
\end{equation*}
which dose not exist.

Above all, we have $\kappa=constant$. Thus we conclude that the conformal transformation must be homothetic.
$\hfill\blacksquare$

%################################################################################

\textbf{Acknowledgement} This work was supported partially by National Natural Science Foundation of China(No.11961061,11761069). The second author would like to thank Professor Bin Chen for his helpful discussion and the valuable comments.\\

\vskip 5mm
\noindent\author{{ \small Cuiling Ma  }\\
{\small {\it College of Mathematics and Systems Science, Xinjiang University}}\\
 {\small {\it Urumqi, Xinjiang Province, 830017, P.R.China}}\\
{\small{\it E-mail: mcl1024@stu.xju.edu.cn}}}

\noindent\author{{ \small Xiaoling Zhang }\\
{\small {\it College of Mathematics and Systems Science, Xinjiang University}}\\
 {\small {\it Urumqi, Xinjiang Province, 830017, P.R.China}}\\
{\small{\it E-mail: zhangxiaoling@xju.edu.cn}}}

\vspace{3mm}


\begin{thebibliography}{99}                                                                                               %
\bibitem  {Antonelli1993} Antonelli P L, Ingarden R, Matsumoto M. The Theory of Sprays and Finsler Space with Applications in Physics and Biology. Kluwer Acad. Publ., Netherlands, 1993.

\bibitem  {Aranasiu2009} Aranasiu G, Neagu M. On Cartan spaces with the m-th root metric $K(x,p)=\sqrt[m]{a^{i_{1} i_{2}...i(m)}(x) p_{i_{1}} p_{i_{2}}\cdot \cdot \cdot p_{i_{m}}}$. Hypercomplex Numbers in Geometry and Physics, 2009, 12:67-73.

\bibitem  {Brinzei2008} Brinzei N. Projective relations for m-th root metric spaces. arXiv:0711.4781, 2008.

%\bibitem  {Chen2013} Chen G, Cheng X. An important class of conformally flat weak Einstein Finsler metrics. Internat. J. Math., 24(2013),no. 1,1350003,15pp.

%\bibitem  {Cheng2014} Cheng X, Li H, Zou Y. On conformally flat $(\alpha,\beta)$-metrics with relatively isotropic mean Landsberg curvature. Publ. Math. Debrecen, 2014, 85:131-144.

\bibitem  {Cheng2012} Cheng X, Shen Z, Tian Y. A class of Einstein $(\alpha,\beta)$-metrics. Israel J. Math., 2012, 192:221-249.

\bibitem  {Chen2015} Chen G, He Q, Shen Z. On conformally flat $(\alpha,\beta)$-metrics with constant flag curvature. Publ. Math. Debrecen, 2015, 86:387-400.

\bibitem {Chen2019} Chen B, Xia K. On conformally flat polynomial $(\alpha,\beta)$-metrics with weakly isotropic scalar curvature. J. Mathematics, 2019, 56:329-352.

\bibitem  {Hashiguchi1976} Hashiguchi M. On conformal transformations of Finsler metrics. J. Math. Kyoto Univ., 1976, 16(1):25-50.

%\bibitem  {Hashuiguchi1977} Hashuiguchi M, Ichijyo Y. On conformal transformations of Wagner spaces. Rep. Fac. Sci. Kagoshima Univ., 1977, 10:19-25.

%\bibitem  {Hashuiguchi1989} Hashuiguchi M, Ichijyo Y. On the condition that a Randers space be conformally flat. Rep. Fac. Sci. Kagoshima Univ. Math. Phys. Chem., 1989, 22:7-14.

%\bibitem  {Hojo2000} Hojo S, Matsumoto M, Okubo K. Theory of conformally Berwald Finsler space and its applications to $(\alpha,\beta)$-metrics. Balkan J. Geom. Appl., 2000, 5:107-118.

\bibitem  {Kang2011} Kang L. On conformally flat Randers metrics. Sci. Sin. Math., 2011, 41:439-446.

%\bibitem  {Kikuchi1994} Kikuchi S. On the condition that a Finsler space be conformally flat. Tensor(N.S), 1994, 55:97-100.

\bibitem  {Knebelman1929} Knebelman M S. Conformal geometry of generalised metric spaces. Proc. Natl. Acad. Sci. USA, 1929, 15:376-379.

%\bibitem  {Li2012} Li B, Shen Z. Projectively flat fourth root Finsler metrics. Canadian Mathematical Bulletin, 2012, 55(1):138-145.

%\bibitem  {Matsumoto1994} Matsumoto M. Conformal change of Finsler space with 1-form metric. An. Stiint. Univ. AI. I. Cuza Iaasi Sect. I a Mar., 1994, 40(1):97-102.

\bibitem  {Matsumoto1996} Matsumoto M. Theory of Finsler spaces with m-th root metric. Publ. Math. Debrecen, 1996, 49:135-155.

%\bibitem  {Matsumoto2001} Matsumoto M. Conformally Berwald and conformally flat Finsler spaces. Publ. Math. Debrecen. 2001, 58:275-285.

\bibitem  {Matsumoto1979} Matsumoto M, Numara S. On Finsler spaces with a cubic metric. Tensor(N.S.), 1979, 38:153-162.

\bibitem  {Rund1959} Rund H. The Differential Geometry of Finsler Spaces. Springer-Verlag, Berlin, 1959.

\bibitem  {Shimada1979} Shimada H. On Finsler spaces with metric $L=\sqrt[m]{a_{i_{1}i_{2...}i_{m}} y^{i_{1}} y^{i_{2}}\dots y^{i_{m}}}$. Tensor, N. S., 1979, 33:365-372.

%\bibitem  {Tayebi2018} Tayebi A, Razgordani M. Four families of projectively flat Finsler metrics with K=1 and their non-Riemannian curvature properties. Rev. R. Acad. Cienc. Exactas FIS. Nat. Ser. A Math. RACSAM., 2018, 112:1463-1485.

\bibitem {Tayebi2019} Tayebi A, Razgordani  M, Najafi B. On conformally flat cubic $(\alpha,\beta)$-metrics. Differ. Geom. Appl., 2019, 62:1-22.

\bibitem  {Wegener1935} Wegener J M. Untersuchungen der zwei- und dreidimensionalen Finslerschen Raumemit der Grundform $L=\sqrt[3]{a_{ijk} x^{'i} x^{'k} x^{'l}}$. Koninkl. Akad. Wetensch, Amsterdam, Proc., 1935, 38:949-955.

\bibitem {Yu2010} Yu Y, You Y. On Einstein m-th root metrics. Differ. Geom. Appl., 2010, 28:290-294.

\bibitem {Zhang2014} Zhang X, Xia Q. On Einstein Matsumoto metrics. Sci. China. Math., 2014, 57:1517-1524.
\end{thebibliography}
\end{document}